\documentclass[letter,10pt]{amsart}
\usepackage{latexsym}
\usepackage{graphics}
\usepackage{amssymb}
\usepackage{epsfig}

\begin{document}
\def\vt{\vartriangle}
\def\bn{$(Bin(X),$ $ \Box)$ }
\def\bn{$(Bin(X),$ $ \Box)$ }
\def\cvb{
\begin{tabular}{c | c c }
      $\bullet$ & x & y   \\
  \hline
   x & x & y \\
   y & a & y   \\
   \end{tabular}
}
\def\ccb{
\begin{tabular}{c | c c }
      $*$ & x & y   \\
  \hline
   x & x & x \\
   y & x & y   \\
   \end{tabular}
}
\def\vb{
\begin{tabular}{c | c c }
      $\bullet$ & x & y   \\
  \hline
   x & x & y \\
   y & x & y   \\
   \end{tabular}
}
\def\BMH{
\begin{tabular}{c|c c c}
      $\bullet$ & a & b & c  \\
  \hline
   a & a & b & c \\
   b & b & b & c  \\
   c & c & b & c \\
  \end{tabular}
}
\def\Tus{
\begin{tabular}{c | c c }
      $\bullet$ & a & b   \\
  \hline
   a & a & b \\
   b & b & b   \\
   \end{tabular}
}
\def\cal{
\begin{tabular}{c | c c }
      $\bullet$ & b & c   \\
  \hline
   b & b & c \\
   c & b & c   \\
   \end{tabular}
}
\def\oosa{
\begin{tabular}{c | c c }
      $\bullet$ & a & c   \\
  \hline
   a & a & c \\
   c & c & c   \\
   \end{tabular}
}
\def\as{
\begin{tabular}{c|c c c}
      $\bullet$ & a & b & c  \\
  \hline
   a & a & b & c \\
   b & b & b & c  \\
   c & c & b & c \\
  \end{tabular}
}
\def\BM2{
\begin{tabular}{c|c c c}
      $\bullet$ & a & b & c  \\
  \hline
   a & a & b & c \\
   b & b & b & b  \\
   c & c & c & c \\
  \end{tabular}
}
\def\bg{
\begin{tabular}{c|c c c}
      $\bullet$ & a & b & c  \\
  \hline
   a & a & b & c \\
   b & b & b & c  \\
   c & c & a & c \\
  \end{tabular}
}
\def\vf{
\begin{tabular}{c|c c c}
      $\bullet$ & a & b & c  \\
  \hline
   a & a & b & c \\
   b & b & b & a  \\
   c & c & c & c \\
  \end{tabular}
}
\def\cd{
\begin{tabular}{c|c c c}
      $\bullet$ & a & b & c  \\
  \hline
   a & a & b & c \\
   b & c & b & b  \\
   c & c & c & c \\
  \end{tabular}
}
\def\xs{
\begin{tabular}{c|c c c}
      $\bullet$ & a & b & c  \\
  \hline
   a & a & c & c \\
   b & b & b & b  \\
   c & c & c & c \\
  \end{tabular}
}
\def\ta{
\begin{tabular}{c | c c }
      $\Box$ & L & R   \\
  \hline
   L & L & R \\
   R & R & L   \\
   \end{tabular}
}

\title{Locally-zero Groupoids and the Center of $Bin(X)$ }
\author{Hiba F. Fayoumi}
\maketitle

\begin{abstract}
In this paper we introduce the notion of the center $ZBin(X)$ in the semigroup
$Bin(X)$ of all binary systems on a set $X$, and show that if $(X,\bullet
)\in ZBin(X)$, then $x\not=y$ implies $\{x,y\}=\{x\bullet y,y\bullet x\}$.
Moreover, we show that a groupoid $(X,\bullet )\in ZBin(X)$ if and only if
it is a locally-zero groupoid.
\end{abstract}

\vspace{10mm}

\bigskip

\markboth{\footnotesize \rm  Hiba F. Fayoumi }
 {\footnotesize \rm Locally-zero Groupoids and the Center of $Bin(X)$  }

\bigskip

\bigskip

\renewcommand{\thefootnote}{}

\medskip\medskip 

\section{\textbf{Preliminaries} \protect\bigskip}

The notion of the semigroup $(Bin(X),$ $\Box )$ was introduced by H. S. Kim
and J. Neggers \textrm{([4])}. Given binary operations \textquotedblleft $%
\ast $" and \textquotedblleft $\bullet $" on a set $X$, they defined a
product binary operation \textquotedblleft $\Box $" as follows: $x\Box
y:=(x\ast y)\bullet (y\ast x)$. This in turn yields a binary operation on $%
Bin(X)$, the set of all groupoids defined on $X$ turning $(Bin(X),\Box )$
into a semigroup with identity ($x\ast y=x$), the left-zero-semigroup, and
an analog of negative one in the right-zero-semigroup. \medskip

\textbf{Theorem 1.1\textrm{{([2])}.}} \textsl{The collection $(Bin(X),$ $%
\Box )$ of all binary systems (groupoids or algebras) defined on $X$ is a
semigroup, i.e., the operation $\Box $ as defined in general is associative.
Furthermore, the left-zero-semigroup is an identity for this operation.}
\medskip

\textbf{Example 1.2\textrm{{([2])}.}} Let $(R,+,\cdot ,0,1)$ be a
commutative ring with identity and let $L(R)$ denote the collection of
groupoids $(R,\ast )$ such that for all $x,y\in R$
\begin{equation*}
x\ast y=ax+by+c
\end{equation*}%
where $a,b,c\in R$ are fixed constants. We shall consider such groupoids to
be \textit{linear groupoids}. Notice that $a=1,b=c=0$ yields $x\ast y=1\cdot
x=x$, and thus the left-zero-semigroup on $R$\ is a linear groupoid. Now,
suppose that $(R,\ast )$ and $(R,\bullet )$ are linear groupoids where $%
x\ast y=ax+by+c$ and $x\bullet y=dx+ey+f$. Then $x\,\Box
\,y=d(ax+by+c)+e(ay+bx=c)+f=(da+eb)x+(db+ea)y+(d+e)c+f$, whence $(R,\Box
)=(R,\ast )\Box (R,\bullet )$ is also a linear groupoid, i.e., $(L(R),\Box )$
is a semigroup with identity. \medskip

\textbf{Example 1.3 \textrm{{([2])}.}} Suppose that in $Bin(X)$ we
consider
all those groupoids $(X,\ast )$ with the \textit{orientation property}: $%
x\ast y\in \{x,y\}$ for all $x$ and $y$. Thus, $x\ast x=x$ as a consequence.
If $(X,\ast )$ and $(X,\bullet )$ both have the orientation property, then
for $x\,\Box \,y=(x\ast y)\bullet (y\ast x)$ we have the possibilities: $%
x\ast x=x,\,y\ast y=y,\,x\ast y\in \{x,y\}$ and $y\ast x\in \{x,y\}$, so
that $x\,\Box \,y\in \{x,y\}$. It follows that if $OP(X)$ denotes this
collection of groupoids, then $(OP(X),\Box )$ is a subsemigroup of $(Bin(X),$
$\Box )$ . In a sequence of papers Nebesk\'{y} \textrm{([3, 4, 5])} has
sought to associate with graphs $(V,E)$ groupoids $(V,\ast )$ with various
properties and conversely. He defined a \textit{travel groupoid} $(X,\ast )$
as a groupoid satisfying the axioms: $(u\ast v)\ast u=u$ and $(u\ast v)\ast
v=u$ implies $u=v$. If one adds these two laws to the orientation property,
then $(X,\ast )$ is an OP-travel-groupoid. In this case $u\ast v=v$ implies $%
v\ast u=u$, i.e., $uv\in E$ implies $vu\in E$, i.e., the digraph $(X,E)$ is
a (simple) graph if $uu\not\in E$, with $u\ast u=u$. Also, if $u\not=v$,
then $u\ast v=u$ implies $(u\ast v)\ast v=u\ast v=u$ is impossible, whence $%
u\ast v=v$ and $uv\in E$, so that $(X,E)$ is a complete (simple) graph.
\bigskip

\section{\textbf{The Center of }$Bin\left( X\right) $\textbf{\ }\protect%
\bigskip}

Let $(X,*)$ be a groupoid and let $ZBin(X)$ denote the collection of
elements of $Bin(X)$ such that $(X,*)\,\Box\,(X,\bullet)\,= \,
(X,\bullet)\,\Box\, (X,*),$ $\forall (X,*)\in Bin(X)$. We call $ZBin(X)$ a
\textit{center} of the semigroup $Bin(X)$. \medskip

\noindent \textbf{Proposition 2.1.} \textsl{The left-zero-semigroup and the
right-zero-semigroup on $X$ are both in $ZBin(X)$.} \medskip

\noindent \textit{Proof.} Given a groupoid $(X,\ast )$, let $(X,\bullet )$
be a left-zero-semigroup. Then $(x\bullet y)\ast (y\bullet x)=x\ast y=(x\ast
y)\bullet (y\ast x)$ for all $x,y\in X$, proving $(X,\ast )\in ZBin(X)$.
Similarly, it holds for the right-zero-semigroup. $\blacksquare $\medskip

\noindent \textbf{Proposition 2.2.} \textsl{If $(X,\bullet )\in ZBin(X)$,
then $x\bullet x=x$ for all $x\in X$.} \medskip

\noindent \textit{Proof.} If $(X,\bullet )\in ZBin(X)$, then $(X,\bullet
)\Box (X,\ast )=(X,\ast )\Box (X,\bullet )$ for all $(X,\ast )\in Bin(X)$.
Let $(X,\ast )\in Bin(X)$ defined by $x\ast y=a$ for any $x,y\in X$ where $%
a\in X$. Then $(x\bullet y)\ast (y\bullet x)=a$ and $(x\ast y)\bullet (y\ast
x)=a\bullet a$ for any $x,y\in X$. Hence we obtain $a\bullet a=a$. If we
change $\left( X,\ast \right) $ in $Bin\left( X\right) $ so that $x\ast y=b$
for every $x,y\in X$ and $b$ is any other element of $X$, then we find that $%
a\bullet a=a$ for any $a\in X$. $\blacksquare $\medskip

Any set can be well-ordered by well-ordering principle, and a well-ordered
set is linearly ordered. With this notion we prove the following. \medskip

\noindent \textbf{Theorem 2.3.} \textsl{If $(X,\bullet )\in ZBin(X)$, then $%
x\not=y$ implies $\{x,y\}=\{x\bullet y,y\bullet x\}$} \medskip

\noindent \textit{Proof.} Let $(X,<)$ be a linear ordered set and let $%
(X,\ast )\in Bin(X)$ be defined by
\begin{equation*}
x\ast y:=\min \{x,y\},\,\,\,\forall x,y\in X\eqno\left( 1\right)
\end{equation*}%
Then we have the following:%
\begin{equation*}
(x\ast y)\bullet (y\ast x)=%
\begin{cases}
x & \text{ if $x\leq y$} \\
y & \text{otherwise}%
\end{cases}%
\eqno\left( 2\right)
\end{equation*}%
\noindent Similarly, we have
\begin{equation*}
(x\bullet y)\ast (y\bullet x)=\min \{x\bullet y,y\bullet x\}\in \{x\bullet
y,y\bullet x\}\eqno\left( 3\right)
\end{equation*}%
If $(X,\bullet )\in ZBin(X)$, then $x<y$ implies $x\in \{x\bullet y,y\bullet
x\}$ for all $x,y\in X$. Similarly, if we define $(X,\ast )\in Bin(X)$ by $%
x\ast y:=\max \{x,y\}$, for all $x,y\in X$, then $x<y$ implies $x\in
\{x\bullet y,y\bullet x\}$ for all $x,y\in X$ when $(X,\bullet )\in ZBin(X)$%
. In any case, we obtain that if $(X,\bullet )\in ZBin(X)$, then
\begin{equation*}
x,y\in \{x\bullet y,y\bullet x\}\eqno\left( 4\right)
\end{equation*}%
We consider 4 cases: (i) $x<y,x\bullet y<y\bullet x$; (ii) $x<y,y\bullet
x<x\bullet y$; (iii) $y<x,x\bullet y<y\bullet x$; (iv) $y<x,y\bullet
x<x\bullet y$. Routine calculations give us the conclusion that $%
\{x,y\}=\{x\bullet y,y\bullet x\}$. $\blacksquare $\medskip

\noindent \textbf{Proposition 2.4.} \textsl{Let $(X,\bullet )\in ZBin(X)$.
If $x\not=y$ in $X$, then $(\{x,y\},\bullet )$ is either a
left-zero-semigroup or a right-zero-semigroup.} \medskip

\noindent \textit{Proof.} Assume that $(X,\bullet )$ is not a
left-zero-semigroup and $x\not=y$ in $X$. Then $(X,\bullet )$ has a subtable:

\begin{equation*}
\begin{tabular}{c|cc}
$\bullet $ & $x$ & $y$ \\ \hline
$x$ & $x$ & $y$ \\
$y$ & $a$ & $y$%
\end{tabular}%
\end{equation*}%
where $a\in \{x,y\}$. Note that $x\bullet x=x,y\bullet y=y$ by Proposition
2.2. Let $(X,\ast )\in Bin(X)$ such that $X$ has a subtable:

\begin{equation*}
\begin{tabular}{c|cc}
$\ast $ & $x$ & $y$ \\ \hline
$x$ & $x$ & $x$ \\
$y$ & $x$ & $y$%
\end{tabular}%
\end{equation*}%
Since $(X,\bullet )\in ZBin(X)$, we have $(x\ast y)\bullet (y\ast
x)=(x\bullet y)\ast (y\bullet x)$ and hence $x\bullet x=y\ast a$. If $a=x$,
then $x\bullet x=y\ast x=x$. If $a=y$, then $x=x\bullet x=y\ast y=y$, a
contradiction. Hence $(X,\bullet )$ should have a subtable:

\begin{equation*}
\begin{tabular}{c|cc}
$\bullet $ & $x$ & $y$ \\ \hline
$x$ & $x$ & $y$ \\
$y$ & $x$ & $y$%
\end{tabular}%
\end{equation*}%
This means $(X,\bullet )$ should be a right-zero-semigroup. Similarly, if $%
(X,\bullet )$ is not a right-zero-semigroup, then it must have a $2\times 2$
table of a left-zero-semigroup. $\blacksquare $\medskip

\noindent \textbf{Proposition 2.5.} \textsl{If $(\{x,y\},\bullet )$ is
either a left-zero-semigroup or a right-zero-semigroup for any $x\not=y$ in $%
X$, then $(X,\bullet )\in ZBin(X)$.} \medskip

\noindent \textit{Proof.} Given $(X,\ast )\in Bin(X)$, let $x\not=y$ in $X$.
Consider $(x\ast y)\bullet (y\ast x)$ and $(x\bullet y)\ast (y\bullet x)$.
If we assume that $(\{x,y\},\bullet )$ is a left-zero-semigroup, then $%
(x\ast y)\bullet (y\ast x)=x\ast y=(x\bullet y)\ast (y\bullet x)$.
Similarly, if we assume that $(\{x,y\},\bullet )$ is a right-zero-semigroup,
then $(x\ast y)\bullet (y\ast x)=y\ast x=(x\bullet y)\ast (y\bullet x)$.
Hence $(X,\bullet )\in ZBin(X)$. $\blacksquare $\medskip

\noindent \textbf{Example 2.6.} Let $X:=\{a,b,c\}$ with the following table:

\begin{equation*}
\begin{tabular}{r|rrr}
$\bullet $ & $a$ & $b$ & $c$ \\ \hline
$a$ & $a$ & $a$ & $c$ \\
$b$ & $b$ & $b$ & $b$ \\
$c$ & $a$ & $c$ & $c$%
\end{tabular}%
\end{equation*}%
Then $(X,\bullet )$ is neither a left-zero-semigroup nor a
right-zero-semigroup, while it has the following subtables:%
\begin{equation*}
\begin{tabular}{r|rr}
$\bullet $ & $a$ & $b$ \\ \hline
$a$ & $a$ & $a$ \\
$b$ & $b$ & $b$%
\end{tabular}%
\quad
\begin{tabular}{r|rr}
$\bullet $ & $a$ & $c$ \\ \hline
$a$ & $a$ & $c$ \\
$c$ & $a$ & $c$%
\end{tabular}%
\quad
\begin{tabular}{r|rr}
$\bullet $ & $b$ & $c$ \\ \hline
$b$ & $b$ & $b$ \\
$c$ & $c$ & $c$%
\end{tabular}%
\end{equation*}

By applying Proposition 2.5, we can see that $(X,\bullet )\in ZBin(X)$.
\medskip

\noindent \textbf{Proposition 2.7.} \textsl{Let }$Ab\left( X\right) $\textsl{%
\ be the collection of all commutative binary systems on }$X$\textsl{. Then }%
$Ab\left( X\right) $\textsl{\ is a right ideal of }$ZBin\left( X\right) $%
\textsl{.}

\noindent \textit{Proof.} Let $\left( X,\bullet \right) \in ZBin\left(
X\right) $ and $\left( X,\ast \right) \in Ab\left( X\right) $. Then by
Proposition 2.2, we have $x\Box y=\left( x\ast y\right) \bullet \left( y\ast
x\right) =\left( x\ast y\right) \bullet \left( x\ast y\right) =x\ast y$.
Also, by Proposition 2.2, we get $y\Box x=\left( y\ast x\right) \bullet
\left( x\ast y\right) =\left( y\ast x\right) \bullet \left( y\ast x\right)
=y\ast x$. Therefore, $\left( X,\ast \right) \Box \left( X,\bullet \right)
\in Ab\left( X\right) $ and so $Ab\left( X\right) \Box ZBin\left( X\right)
\subseteq Ab\left( X\right) $. $\blacksquare $\medskip

\section{\textbf{Locally-zero Groupoids}\protect\bigskip}

A groupoid $(X,\bullet )$ is said to be \textit{locally-zero} if (i) $%
x\bullet x=x$ for all $x\in X$; (ii) for any $x\not=y$ in $X$, $%
(\{x,y\},\bullet )$ is either a left-zero-semigroup or a
right-zero-semigroup. \medskip

Using Propositions 2.2, 2.4 and 2.5 we obtain the following. \medskip

\noindent \textbf{Theorem 3.1.} \textsl{A groupoid $(X,\bullet )\in ZBin(X)$
if and only if it is a locally-zero groupoid.} \medskip

Given any two elements $x,y\in X$, there exists exactly one
left-zero-semigroup and one right-zero-semigroup, and so if we apply Theorem
3.1 we have the following Corollary. \medskip

\noindent \textbf{Corollary 3.2.} \textsl{If $|X|=n$, there are $2^{\binom{n%
}{2}}$-different (but may not be isomorphic) locally-zero groupoids.}
\medskip

\noindent For example, if $n=3$, there are $2^{3}=8$ such groupoids, i.e.,

\begin{equation*}
\begin{tabular}{c|ccc}
$\bullet $ & $a$ & $b$ & $c$ \\ \hline
$a$ & $a$ & $a$ & $a$ \\
$b$ & $b$ & $b$ & $b$ \\
$c$ & $c$ & $c$ & $c$%
\end{tabular}%
\ \quad
\begin{tabular}{c|ccc}
$\bullet $ & $a$ & $b$ & $c$ \\ \hline
$a$ & $a$ & $a$ & $a$ \\
$b$ & $b$ & $b$ & $c$ \\
$c$ & $c$ & $b$ & $c$%
\end{tabular}%
\text{ }\quad
\begin{tabular}{c|ccc}
$\bullet $ & $a$ & $b$ & $c$ \\ \hline
$a$ & $a$ & $a$ & $c$ \\
$b$ & $b$ & $b$ & $b$ \\
$c$ & $a$ & $c$ & $c$%
\end{tabular}%
\text{ }\quad
\begin{tabular}{c|ccc}
$\bullet $ & $a$ & $b$ & $c$ \\ \hline
$a$ & $a$ & $b$ & $a$ \\
$b$ & $a$ & $b$ & $b$ \\
$c$ & $c$ & $c$ & $c$%
\end{tabular}%
\end{equation*}%
\medskip

\begin{equation*}
\begin{tabular}{c|ccc}
$\bullet $ & $a$ & $b$ & $c$ \\ \hline
$a$ & $a$ & $b$ & $c$ \\
$b$ & $a$ & $b$ & $c$ \\
$c$ & $a$ & $b$ & $c$%
\end{tabular}%
\ \quad
\begin{tabular}{c|ccc}
$\bullet $ & $a$ & $b$ & $c$ \\ \hline
$a$ & $a$ & $a$ & $c$ \\
$b$ & $b$ & $b$ & $c$ \\
$c$ & $a$ & $b$ & $c$%
\end{tabular}%
\ \quad
\begin{tabular}{c|ccc}
$\bullet $ & $a$ & $b$ & $c$ \\ \hline
$a$ & $a$ & $b$ & $a$ \\
$b$ & $a$ & $b$ & $c$ \\
$c$ & $c$ & $b$ & $c$%
\end{tabular}%
\ \quad
\begin{tabular}{c|ccc}
$\bullet $ & $a$ & $b$ & $c$ \\ \hline
$a$ & $a$ & $b$ & $c$ \\
$b$ & $a$ & $b$ & $b$ \\
$c$ & $a$ & $c$ & $c$%
\end{tabular}%
\end{equation*}%
\medskip

\noindent \textbf{Corollary 3.3.} \textsl{The collection of all locally-zero
groupoids on $X$ forms a subsemigroup of $(Bin(X),\Box )$.} \medskip

\noindent \textit{Proof.} Let $x\not=y$ in $X$. If $(\{x,y\},\bullet )$ is a
left-zero-semigroup and $(\{x,y\},\ast )$ is a right-zero semigroup, then $%
x\Box y=(x\bullet y)\ast (y\bullet x)=x\ast y=y$, $y\Box x=(y\bullet x)\ast
(x\bullet y)=y\ast x=x$, i.e., $(\{x,y\},\Box )$ is a right-zero-semigroup.
Similarly, we can prove the other three cases, i.e.,

\begin{equation*}
\begin{tabular}{c|cc}
$\Box $ & $L$ & $R$ \\ \hline
$L$ & $L$ & $R$ \\
$R$ & $R$ & $L$%
\end{tabular}%
\end{equation*}%
where $L$ means the \textquotedblleft left-zero-semigroup" and $R$ means
that \textquotedblleft right-zero-semigroup", proving that the locally-zero
groupoids on $X$ form a subsemigroup of $(Bin(X),\Box )$. $\blacksquare $%
\medskip

Using Corollary 3.3, we can see that $(X,\bullet )\Box (X,\ast )$ belongs to
the center $ZBin\left( X\right) $ of $Bin(X)$ for any $(X,\bullet ),(X,\ast
)\in ZBin(X)$. \medskip

\noindent \textbf{Proposition 3.4. }\textsl{Not all locally-zero groupoids
are semigroups.}\medskip

\noindent \textit{Proof. }Consider $\left( X,\bullet \right) $ where $%
X:=\left\{ a,b,c\right\} $ and $"\bullet "$ is given by the following table:
\begin{equation*}
\begin{tabular}{r|rrr}
$\bullet $ & $a$ & $b$ & $c$ \\ \hline
$a$ & $a$ & $a$ & $c$ \\
$b$ & $b$ & $b$ & $b$ \\
$c$ & $a$ & $c$ & $c$%
\end{tabular}%
\end{equation*}%
Then it is easy to see that $\left( X,\bullet \right) $ is locally-zero.
Consider the subtables:
\begin{equation*}
\begin{tabular}{r|rr}
$\bullet $ & $a$ & $b$ \\ \hline
$a$ & $a$ & $a$ \\
$b$ & $b$ & $b$%
\end{tabular}%
\quad \text{ }%
\begin{tabular}{r|rr}
$\bullet $ & $a$ & $c$ \\ \hline
$a$ & $a$ & $c$ \\
$c$ & $a$ & $c$%
\end{tabular}%
\quad \text{ }%
\begin{tabular}{r|rr}
$\bullet $ & $b$ & $c$ \\ \hline
$b$ & $b$ & $b$ \\
$c$ & $c$ & $c$%
\end{tabular}%
\end{equation*}%
and notice that $\left( \left\{ a,b\right\} ,\bullet \right) $, $\left(
\left\{ a,c\right\} ,\bullet \right) $ and $\left( \left\{ b,c\right\}
,\bullet \right) $ are left-, right- and left-zero-semigroups, respectively.
But $\left( a\bullet b\right) \bullet c=a\bullet c=c,$ while $a\bullet
\left( b\bullet c\right) =a\bullet b=a$. Hence $\left( X,\bullet \right) $
fails to be a semigroup and the result follows. $\blacksquare $\medskip

\noindent \textbf{Proposition 3.5}. \textsl{Let }$\left( X,\bullet \right) $%
\textsl{\ be a locally-zero groupoid. If }$\left( X,\bullet \right) $\textsl{%
\ is a semigroup then it is either a left- or a right-zero-semigroup.}%
\medskip

\noindent \textit{Proof.} Suppose that $\left( X,\bullet \right) $ is a
semigroup, then $\left( x\bullet y\right) \bullet z=x\bullet \left( y\bullet
z\right) $ for all $x,y,z\in X$. By Theorem 3.1, $\left( X,\bullet \right) $
$\in ZBin\left( X\right) $, and then by Proposition 2.4, $\left( \left\{
x,y\right\} ,\bullet \right) $ is either a left-zero- or a
right-zero-semigroup for $x\neq y$. In fact, $\left( \left\{ x,z\right\}
,\bullet \right) $ and $\left( \left\{ y,z\right\} ,\bullet \right) $ are
also either left-zero- or right-zero-semigroups for $x\neq z$ and $y\neq z,$
respectively. Assume that $\left( \left\{ x,y\right\} ,\bullet \right) $, $%
\left( \left\{ x,z\right\} ,\bullet \right) $ and $\left( \left\{
y,z\right\} ,\bullet \right) $ are left-, right- and left-zero-semigroups,
respectively. Then, $\left( x\bullet y\right) \bullet z=x\bullet z=z$ while $%
x\bullet \left( y\bullet z\right) =x\bullet y=x$, a contradiction.
Similarly, we can reach a contradiction if we assume that $\left( \left\{
x,y\right\} ,\bullet \right) $, $\left( \left\{ x,z\right\} ,\bullet \right)
$ and $\left( \left\{ y,z\right\} ,\bullet \right) $ are right-, left- and
right-zero-semigroups, respectively. Now suppose that $\left( \left\{
x,y\right\} ,\bullet \right) $, $\left( \left\{ x,z\right\} ,\bullet \right)
$ and $\left( \left\{ y,z\right\} ,\bullet \right) $ are left-, left- and
right-zero-semigroups, respectively. Then, $\left( y\bullet x\right) \bullet
z=y\bullet z=z$ while $y\bullet \left( x\bullet z\right) =y\bullet x=y$, a
contradiction. Similarly, we can reach a contradiction if we assume that $%
\left( \left\{ x,y\right\} ,\bullet \right) $, $\left( \left\{ x,z\right\}
,\bullet \right) $ and $\left( \left\{ y,z\right\} ,\bullet \right) $ are
right-, right- and left-zero-semigroups, respectively. Hence, the only two
other cases are when all three subgroupoids are either all left- or all
right-zero-semigroups. Therefore, $\left( X,\bullet \right) $ is either
left- or right-zero-semigroup. $\blacksquare $\medskip

\noindent \textbf{Proposition 3.6. }\textsl{Let }$\left( X,\bullet \right) $%
\textsl{\ be a locally-zero groupoid. Then }$\left( X,\bullet \right) \Box
\left( X,\bullet \right) =\left( X,\Box \right) $\textsl{\ is the
left-zero-semigroup on }$X$\textsl{.}\medskip

\noindent \textit{Proof. }Suppose that $\left( \left\{ x,y\right\} ,\bullet
\right) $ is the right-zero-semigroup, then $x\Box y=\left( x\bullet
y\right) \bullet \left( y\bullet x\right) =y\bullet x=x.$On the other hand,
if $\left( \left\{ x,y\right\} ,\bullet \right) $ is the
left-zero-semigroup, then $x\Box y=\left( x\bullet y\right) \bullet \left(
y\bullet x\right) =x\bullet y=x.$Thus in both cases, $x\Box y=x$ for all $%
x\in X$ making $\left( X,\Box \right) $ the left-zero-semigroup. $%
\blacksquare $

\bigskip\bigskip

{\footnotesize {\ \textsc{\ Hiba Fayoumi, Department of Mathematics,
University of Alabama, Tuscaloosa, AL, 35487-0350, U. S. A.} }}

{\footnotesize \textit{E-mail address}: \texttt{hiba.fayoumi@ua.edu} }

\end{document}